\documentclass{gtpart}     
\usepackage{pinlabel}  
\usepackage{graphicx}  

\usepackage{amsmath}
\usepackage{amssymb}
\usepackage{amsthm}
\usepackage{tabularx}
\usepackage{longtable}
\usepackage{booktabs}
\usepackage{subcaption} 

\title{From Heegaard splittings to trisections; porting $3$-dimensional ideas to dimension $4$}

\author{David T Gay}
\givenname{David T}
\surname{Gay}
\address{Euclid Lab\\ 160 Milledge Terrace\\ Athens, GA 30606\\\newline Department of Mathematics\\ University
  of Georgia\\ Athens, GA 30602}
\email{d.gay@euclidlab.org}

\volumenumber{}
\issuenumber{}
\publicationyear{}
\papernumber{}
\startpage{}
\endpage{}
\doi{}
\MR{}
\Zbl{}
\received{}
\revised{}
\accepted{}
\published{}
\publishedonline{}
\proposed{}
\seconded{}
\corresponding{}
\editor{}
\version{}

\newtheorem{theorem}{Theorem}
\newtheorem{lemma}[theorem]{Lemma}
\newtheorem{proposition}[theorem]{Proposition}

\theoremstyle{definition}
\newtheorem{definition}[theorem]{Definition}

\def\Z{\mathbb Z}

\def\R{\mathbb R}
\def\C{\mathbb C}

\newcommand{\into}{\ensuremath{\hookrightarrow}}

\newcommand{\Ss}{\mathcal S}
\newcommand{\Tt}{\mathcal T}
\newcommand{\Dd}{\mathcal D}

\begin{document}

\begin{abstract}    
 These notes summarize and expand on a mini-course given at CIRM in February 2018 as part of Winter Braids VIII. We somewhat obsessively develop the slogan ``Trisections are to $4$--manifolds as Heegaard splittings are to $3$--manifolds'', focusing on and clarifying the distinction between three ways of thinking of things: the basic definitions as decompositions of manifolds, the Morse theoretic perspective and descriptions in terms of diagrams. We also lay out these themes in two important relative settings: $4$--manifolds with boundary and $4$--manifolds with embedded $2$--dimensional submanifolds.
\end{abstract}

\maketitle

\section{Introduction}

All manifolds are smooth in this paper, except that a very mild form of manifold with boundary and corners appears without comment at various places, and the appropriate rounding of corners is assumed without comment.

Most of the content of this paper is in the form of definitions and statements of basic results, and some discussion. There are no proofs; either proofs are suggested as exercises, sometimes with hints, or external references are given. We necessarily present a very limited range of material and hope that this a useful launching point for more in-depth reading and, especially, for new and original research.

At the risk of overdoing it, we maintain a format throughout which heavily emphasizes the parallels between the $3$-- and $4$--dimensional settings. In particular, we use a $2$--column format for most definitions and theorems, with parallel bulleted items for the $4$--dimensional setting on the left and the $3$--dimensional setting on the right; sometimes there is an extra condition in dimension four which does not have a three dimensional analog, in which case to avoid excessive white space we drop the $2$--column format for this last condition. This format is based on the approach taken on the blackboard in the original mini-course, and we hope the experiment is equally effective in printed form.

One goal of these notes is to emphasize the Morse theoretic perspective where it often gets conveniently ignored in other presentations. In principle one can understand everything one needs to know about trisections without thinking Morse theoretically, but this seems to miss an essential piece of the intuition. For this reason, in section~\ref{S:Decomps}, we quickly cover the basic definitions of Heegaard splittings and trisections as decompositions of manifolds somewhat drily and minimally so as to get on to the Morse theory of Section~\ref{S:Morse} quickly.

\section{The basic definitions: decompositions} \label{S:Decomps}

We use the symbol $\#$ for connected sum and $\natural$ for boundary connected sum, so that $\partial(A \natural B) = (\partial A) \# (\partial B)$. Then $\#^n A$ is the connected sum of $n$ copies of $A$, with $\#^0 A = S^m$, when $A$ is a manifold of dimension $m$. Similarly $\natural^n B$ is the boundary connected sum of $n$ copies of $B$, with $\natural^0 B = B^m$, when $B$ is a manifold with connected boundary of dimension $m$. With this in mind we name the following standard manifolds of dimensions $2$, $3$ and $4$:
\begin{itemize}
 \item The standard genus $g$ surface is $\Sigma_g = \#^g (S^1 \times S^1)$.
 \item The standard genus $g$ handlebody is $H_g = \natural^g (S^1 \times B^2)$, with $\partial H_g = \Sigma_g$.
 \item The standard $4$--dimensional $1$--handlebody (of ``genus $k$'') is $Z_k = \natural^k (S^1 \times B^3)$.
\end{itemize}

%

\newlength{\halfdeflength}
\newlength{\wholedeflength}

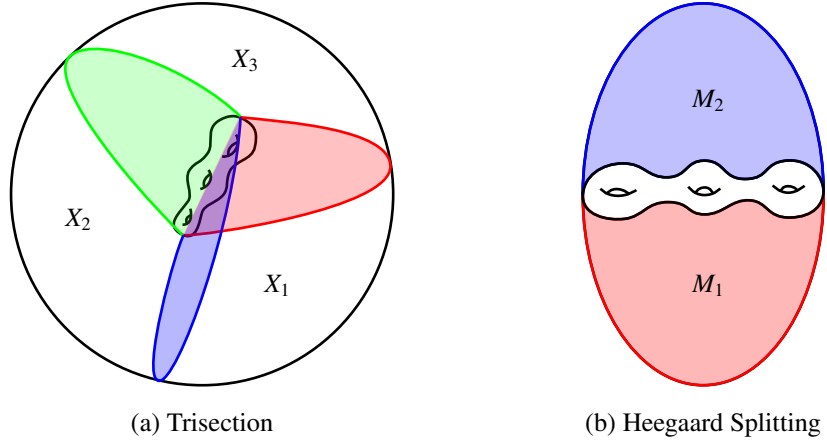
\begin{figure}
  \begin{subfigure}[c]{0.5\textwidth}
    \labellist
    \small\hair 2pt
    \pinlabel $X_1$ at 103 40
    \pinlabel $X_2$ at 28 65
    \pinlabel $X_3$ at 91 125
    \endlabellist
    \centering
    \includegraphics{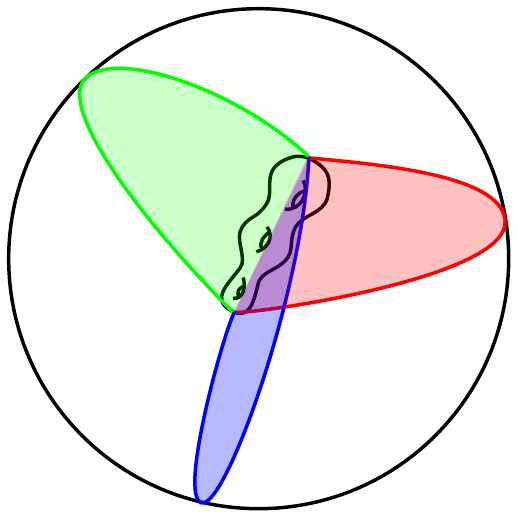}
    \caption{Trisection}
    \label{F:TriSchematic}
  \end{subfigure}
  \begin{subfigure}[c]{0.5\textwidth}
    \labellist
    \small\hair 2pt
    \pinlabel $M_1$ at 50 40
    \pinlabel $M_2$ at 50 110
    \endlabellist
    \centering
    \includegraphics{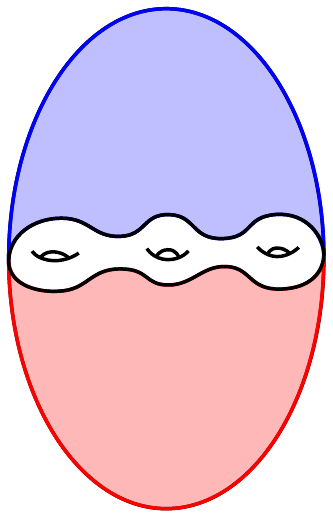}
    \caption{Heegaard Splitting}
    \label{F:HSSchematic}
  \end{subfigure}
  \caption{Schematics of trisections and Heegaard splittings}
   \label{F:HSandTriSchematic}
\end{figure}

\begin{definition} \label{D:HSandTri}
In which we define Heegaard splittings and trisections and establish orientation conventions. (See Figure~\ref{F:HSandTriSchematic}.)
\setlength{\halfdeflength}{0.5\linewidth}
\addtolength{\halfdeflength}{-2\tabcolsep}
\setlength{\wholedeflength}{2\halfdeflength}
\addtolength{\wholedeflength}{2\tabcolsep}
\begin{longtable}{p{\halfdeflength} | p{\halfdeflength}}
\textbf{Dimension four:} A $(g;k_1,k_2,k_3)$ {\em trisection} of a closed, connected, oriented $4$--manifold $X$ is a decomposition $X = X_1 \cup X_2 \cup X_3$ such that: 
&
\textbf{Dimension three:} A genus $g$ {\em Heegaard splitting} of a closed, connected, oriented $3$--manifold $M$ is a decomposition $M = M_1 \cup M_2$ such that: 
\\
\textbullet\; For each $i$, $X_i$ is diffeomorphic to $Z_{k_i}$. 
&
\textbullet\; For each $i$, $M_i$ is diffeomorphic to $H_g$. 
\\
\textbullet\; Taking indices mod $3$, each $X_i \cap X_{i+1}$ is diffeomorphic to $H_g$. We orient $X_i \cap X_{i+1}$ as a submanifold of $\partial X_{i+1}$.
& 
\textbullet\; $M_1 \cap M_2$ is diffeomorphic to $\Sigma_g$. We orient $M_1 \cap M_2$ as $\partial M_1 = -\partial M_2$.
\\
\cmidrule{2-2}
\multicolumn{2}{p{\wholedeflength}}{\textbullet\; $X_1 \cap X_2 \cap X_3$ is diffeomorphic to $\Sigma_g$. We orient $X_1 \cap X_2 \cap X_3$ as $\partial (X_1 \cap X_2) = \partial (X_2 \cap X_3) = \partial (X_3 \cap X_1)$.}
\end{longtable}
A $(g;k_1,k_2,k_3)$ trisection is {\em balanced} if $k_1=k_2=k_3=k$, in which case we call it a $(g,k)$ trisection.
\end{definition}

A Heegaard splitting will often be labelled $\Ss$, to refer to the triple $\Ss=(M, M_1,M_2)$, and similarly a trisection will often be labelled $\Tt$, to refer to the $4$--tuple $\Tt=(X,X_1,X_2,X_3)$. Note that the labeling of the pieces matters; $(M,M_1,M_2)$ and $(M,M_2,M_1)$ are different Heegaard splittings of the same underlying oriented $3$--manifold.

To digest the orientation conventions, a good exercise is to verify first that, in a trisection $\Tt=(X,X_1,X_2,X_3)$, the orientations of $\Sigma = X_1 \cap X_2 \cap X_3$ as $\partial (X_1 \cap X_2)$, $\partial (X_2 \cap X_3)$ and $\partial (X_3 \cap X_1)$ really do agree. Then one should verify that this orientation of $\Sigma$ from the $\Tt$ agrees with its orientation as the splitting surface in each of the the Heegaard splittings $\Ss_i = (\partial X_i,X_{i-1} \cap X_i,X_i \cap X_{i+1})$.

\begin{definition} \label{D:Stab}
 In which we define a stabilization operation for both kinds of decompositions.
 \setlength{\halfdeflength}{0.5\linewidth}
 \addtolength{\halfdeflength}{-2\tabcolsep}
 \begin{longtable}{p{\halfdeflength} | p{\halfdeflength}}
  \textbf{Dimension four:} Given a trisection $\Tt=(X,X_1,X_2,X_3)$ of a $4$--manifold $X$ and an index $i \in\Z /3\Z$, an {\em $i$--stabilization} of this trisection is a trisection $\Tt'=(X,X_1',X_2',X_3')$ obtained as follows:
  &
  \textbf{Dimension three:} Given a Heegaard splitting $\mathcal{S}=(M,M_1,M_2)$ of a $3$--manifold $M$, and an index $i \in \Z /2\Z$, an {\em $i$--stabilization} of $\mathcal{S}$ is a Heegaard splitting $\mathcal{S}' = (M,M_1',M_2')$ of $M$ obtained as follows:
  \\
  \textbullet\; Choose an arc $a$ properly embedded and boundary parallel in $X_{i-1} \cap X_{i+1}$, with a regular neighborhood $\nu \cong B^3 \times a$ so that $\nu \cap X_i \cong B^3 \times \partial a$ and $\nu \cap X_1 \cap X_2 \cap X_3 \cong B^2 \times \partial a$.
  &
  \textbullet\; Choose an arc $a$ properly embedded and boundary parallel in $M_{i+1}$, with a regular neighborhood $\nu \cong B^2 \times a$ so that $\nu \cap M_i \cong B^2 \times \partial a$.
  \\
  \textbullet\; Let $X_i' = X_i \cup \nu$.
  &
  \textbullet\; Let $M_i' = M_i \cup \nu$.
  \\
  \textbullet\; Let $X_{i \pm 1}' = X_{i \pm 1} \setminus \mathring{\nu} $.
  &
  \textbullet\; Let $M_{i+1}' = M_{i+1} \setminus \mathring{\nu} $.
 \end{longtable}
\end{definition}

Several comments are in order. The fact that a stabilization of a Heegaard splitting or trisection is again a Heegaard splitting or trisection is a lemma that needs to be proved, and is a worthwhile exercise. In both dimensions, any $i$--stabilization of a given Heegaard splitting or trisection is isotopic to any other. In dimension three, $1$--stabilization and $2$--stabilization are isotopic, and both turn a genus $g$ splitting into a genus $g+1$ splitting. For this reason in dimension three we dispense with the index and simply say ``stabilization'', an operation defined uniquely up to ambient isotopy. In dimension four, any two $i$--stabilizations of the same trisection are isotopic, but an $i$--stabilizations turns a $(g;k_1,k_2,k_3)$ trisection into a $(g+1;k_1',k_2',k_3')$ trisection where $k'_i = k_i+1$ and, for $j \neq i$, $k'_j=k_j$. Thus $1$--stabilization, $2$--stabilization and $3$--stabilization are necessarily different. By stabilization, any trisection can be made balanced, and ``stabilization'' for a balanced trisection means the result of performing one $1$--, one $2$-- and one $3$--stabilization. This balanced stabilization is the stabilization process originally presented in~\cite{GayKirby}.

The basic results (some discussion of their proofs appears in the following section) are:
\begin{theorem}[Existence and Uniqueness] \label{T:ExUn}
 The above decompositions exist and are unique up to stabilization. More precisely:
 \setlength{\halfdeflength}{0.5\linewidth}
 \addtolength{\halfdeflength}{-2\tabcolsep}
 \begin{longtable}{p{\halfdeflength} | p{\halfdeflength}}
  \textbf{Dimension four:}  Every closed, connected, oriented $4$--manifold has a trisection, and any two trisections of the same $4$--manifold become isotopic after some number of stabilizations.~\cite{GayKirby}
  &
  \textbf{Dimension three:} Every closed, connected, oriented $3$--manifold has a Heegaard splitting, and any two Heegaard splittings of the same $3$--manifold become isotopic after some number of stabilizations.\cite{Reidemeister, Singer}
 \end{longtable}
\end{theorem}

\section{The Morse theoretic perspective} \label{S:Morse}

We assume familiarity with basic Morse theory and the connection between Morse functions on manifolds and handle decompositions. We will define some of these basic notions below in certain cases, only for the purpose of establishing parallels between $3$-- and $4$--dimensional phenomena.
\begin{definition} \label{D:MorseMorse2}
In which we define Morse functions and Morse $2$--functions in the limited context of dimensions three and four.
 \setlength{\halfdeflength}{0.5\linewidth}
 \addtolength{\halfdeflength}{-2\tabcolsep}
 \begin{longtable}{p{\halfdeflength} | p{\halfdeflength}}
  \textbf{Dimension four:} A {\em Morse $2$--function} on a $4$--manifold $X$ is a smooth function $f: X \to \R^2$ which, at every point $p \in X$, has one of the following three forms with respect to appropriate local coordinates $(t,x,y,z)$ near $p$ and $(u,v)$ near $f(p)$:
  &
  \textbf{Dimension three:} A {\em Morse function} on a $3$--manifold $M$ is a smooth function $f: M \to \R$ which, at every point $p \in M$, has one of the following two forms with respect to appropriate local coordinates $(x,y,z)$ near $p$ and $u$ near $f(p)$:
  \\
  \textbullet\; $(t,x,y,z) \mapsto (u=t,v=x)$; here $p$ is called a {\em regular point}.
  &
  \textbullet\; $(x,y,z) \mapsto u=x$; here $p$ is called a {\em regular point}.
  \\
  \textbullet\; $(t,x,y,z) \mapsto (u=t, v = \pm x^2 \pm y^2 \pm z^2)$; here $p$ is called a {\em fold point} and $p$ is called {\em definite} or {\em indefinite} according to whether the quadratic form $\pm x^2 \pm y^2 \pm z^2$ is definite or indefinite.
  &
  \textbullet\; $(x,y,z) \mapsto u = \pm x^2 \pm y^2 \pm z^2$; here $p$ is called a {\em critical point}, and the number of $-$'s in the quadratic form $\pm x^2 \pm y^2 \pm z^2$ is the {\em index of $p$}.
  \\
  \cmidrule{2-2}
  \multicolumn{2}{l}{\textbullet\; $(t,x,y,z) \mapsto (u=t, v= x^3 - t x \pm y^2 \pm z^2)$; here $p$ is called a {\em cusp point}.}
 \end{longtable}
 In both cases, a point $q$ in the codomain of $f$ is called a {\em regular value} if all points $p \in f^{-1}(q)$ are regular points, otherwise $q$ is a {\em critical value}. In dimension four, both fold and cusp points are {\em critical points}.
\end{definition}
A good way to think about the connection between Morse functions and Morse $2$--functions is that, locally, a Morse $2$--function looks like time crossed with a generic homotopy betweeen Morse functions. Along a fold we can parametrize things so that we see a single Morse critical point not moving in time, while a cusp corresponds to a birth or death of a cancelling pair of critical points.

\begin{figure}
    \labellist
    \small\hair 2pt
    \pinlabel $0$ [r] at 64 80
    \pinlabel $2$ [r] at 130 125
    \pinlabel $1$ [bl] at 169 150
    \endlabellist
 \centering
 \includegraphics{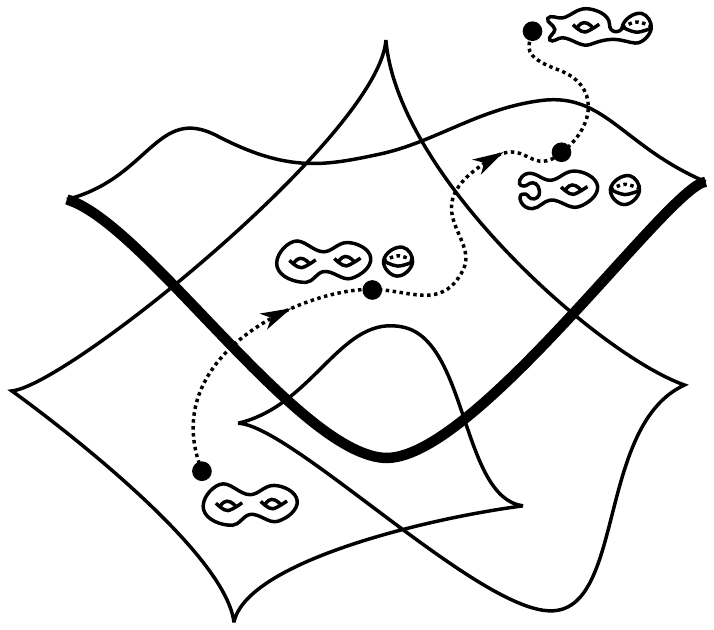}
 \caption{Some characteristic features of a Morse $2$--function. The darker arc is a definite fold and the remaining solid arcs are indefinite folds. The dotted arc is an oriented arc $A$ transverse to the folds, with indices of the critical points of the associated Morse function on $f^{-1}(A)$ labelled at the crossings. The surfaces are representative inverse images of regular values along the arc. \label{F:Morse2FCartoon}}
\end{figure}

Here we recommend that the reader verify the following facts as an exercise in building the correct intution (assume here that the domain of $f$ is closed):
\begin{itemize}
 \item In both cases the inverse image of a regular value is a closed surface.
 \item In both cases the singular locus, the set of all critical points, is a closed codimension three submanifold, i.e. a finite collection of points in dimension three and a finite collection of embedded circles in dimension four.
 \item In dimension four, the cusp points form a finite collection of points on the singular locus.
 \item Via a small perturbation, in dimension three one may assume that the critical points of a Morse function have distinct critical values
 \item Returning to dimension four, letting $Z$ be the singular locus of a Morse $2$--function $f$, via a small perturbation one may assume that $f|_Z$ is an immersion with semicubical cusps, with at worst double point self intersections, none of which occur at cusps. (Figure~\ref{F:Morse2FCartoon} is an attempt at a cartoon illustrating many of the features of a Morse $2$--function discussed in this and the following bullet points.)
 \item If $f:X \to \R^2$ is a Morse $2$--function and $A$ is an arc in $\R^2$ avoiding the cusps and transverse to the image of the singular locus, then $M = f^{-1}(A)$ is a $3$--manifold in $X$, with $\partial M = f^{-1}(\partial A)$. 
 \item Furthermore, if we identify $A$ with an interval in $\R$ via some embedding $A \hookrightarrow \R$ then $f|_M :M \to A$ is a Morse function with critical points of index $0$ and $3$ where $A$ crosses definite folds and critical points of index $1$ and $2$ where $A$ crosses indefinite folds. Reversing the orientation of $A$ changes the indices of these critical points, with index $0$ becoming index $3$ and vice versa, and index $1$ becoming index $2$ and vice versa. 
 \item Crossing a definite fold in the index $0$ direction adds a new $S^2$ component to the fiber (preimage of regular value) while crossing in the index $3$ direction caps off such a component. Crossing an indefinite fold in the index $1$ direction either increases the genus of a fiber component by one or connects two disconnected components, while crossing in the index $2$ direction surgers the fiber along a compressing circle, either decreasing genus by one or splitting a component in two.
 \item If the arc $A$ passes immediately adjacent to a cusp, thus crossing two folds, then the corresponding Morse critical points are two cancelling critical points of successive index.
\end{itemize}

Next we will define the kind of Morse functions and Morse $2$--functions which produce, respectively, Heegaard splittings and trisections. Our $3$--dimensional definition will, as usual, seem a little odd since it is set up to emphasize the parallel with the $4$--dimensional setting. First we introduce some notation. For $\theta \in S^1 = \R /\!2\pi\Z$, let $R_\theta \subset \R^2$ be the ray making angle $-\theta$ with the positive $x$--axis. (Yes, the negative sign in $-\theta$ is intentional, we explicitly want to move {\em clockwise} around the origin, simply because this fits well with other orientation conventions.) Identify $R_\theta$ with $[0,\infty)$ via the parametrization $(t\cos(-\theta), t\sin(-\theta))$. Thinking of $\R$ as the $x$--axis in $\R^2$, the intervals $[0,\infty) \subset \R$ and $(-\infty,0] \subset \R$ are then identified with the rays $R_0$ and $R_\pi$, respectively, except that when we think of $(-\infty,0]$ as $R_\pi$ and then use the above parametrization to identify $R_\pi$ with $[0,\infty)$, we then see $(-\infty,0]$ oriented away from $0$. We also consider the ``trisection'' of $\R^2$ as $\R^2 = A_1 \cup A_2 \cup A_3$ where $A_i$ is the sector bounded by $\R_{2\pi (i-1)/3}$ and $\R_{2\pi i/3}$.
\begin{definition} \label{D:HSTfunctions}
In which we define Heegaard splitting Morse functions and trisecting Morse $2$--functions.
\setlength{\halfdeflength}{0.5\linewidth}
\addtolength{\halfdeflength}{-2\tabcolsep}
\setlength{\wholedeflength}{2\halfdeflength}
\addtolength{\wholedeflength}{2\tabcolsep}
 \begin{longtable}{p{\halfdeflength} | p{\halfdeflength}}
  \textbf{Dimension four:} A {\em $(g;k_1,k_2,k_3)$ trisecting} Morse $2$--function $f$ on a $4$--manifold $X$ is a Morse $2$--function $f: X \to \R^2$ such that:
  &
  \textbf{Dimension three:} A {\em genus $g$ Heegard splitting} Morse function $f$ on a $3$--manifold $M$ is a Morse function $f:M \to \R$ such that:
  \\
  \textbullet\; $\mathbf{0}=(0,0)$ is a regular value of $f$, and thus $f^{-1}(\mathbf{0}) = \Sigma$ is a closed surface, which we require to be connected of genus $g$.
  &
  \textbullet\; $0$ is a regular value of $f$, and thus $f^{-1}(0) = \Sigma$ is a closed surface, which we require to be connected of genus $g$.
  \\
  \textbullet\; On each of the three rays $R_0$, $R_{2\pi/3}$ and $R_{4\pi/3}$, $f$ has exactly $g$ index $2$ and one index $3$ critical points, all of which have distinct critical values.
  &
  \textbullet\; On each of the two rays $R_0$ and $R_\pi$, $f$ has exactly $g$ index $2$ and one index $3$ critical points, all of which have distinct critical values.
  \\
  \cmidrule{2-2}
  \multicolumn{2}{p{\wholedeflength}}{\textbullet\; Over each of sector $A_i$, the singular locus of $f$ has exactly $g+1$ components, all of which are arcs from one bounding ray of $A_i$ to the next, $k_i$ of which are indefinite folds without cusps, $(g-k_i)$ of which are indefinite folds each with exactly one indefinite cusp, and one of which (the outermost) is a definite fold. Furthermore, in $\R^2$ each of these components is transverse to each ray $R_\theta$ except at cusps, which are tangent to the rays, and $f$ restricted to the singular locus is an immersion with cusps and double points avoiding the cusps. This is illustrated in Figure~\ref{F:TrisectingM2F}.}
 \end{longtable}
 
\end{definition}

\begin{figure}
    \labellist
    \small\hair 2pt
    \pinlabel $R_0$ [l] at 157 75
    \pinlabel $R_{2\pi/3}$ [r] at 32 4
    \pinlabel $R_{4\pi/3}$ [r] at 32 146
    \pinlabel $A_1$ at 127 13
    \pinlabel $A_2$ [r] at 0 73
    \pinlabel $A_3$ at 126 138
    \endlabellist
 \centering
 \includegraphics{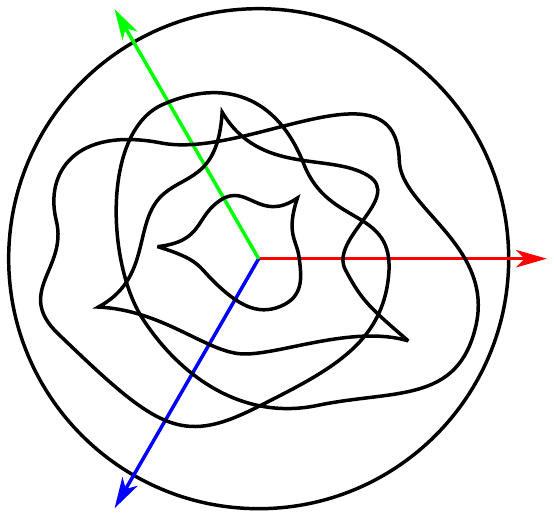}
 \caption{The singular value set of a $(4;3,2,2)$ trisecting Morse $2$--function} \label{F:TrisectingM2F}
\end{figure}
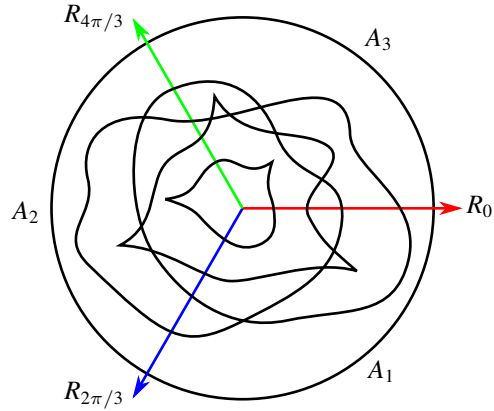

\begin{lemma} \label{L:Fcn2Decomp}
 These types of functions induce the indicated manifold decompositions:
\setlength{\halfdeflength}{0.5\linewidth}
 \addtolength{\halfdeflength}{-2\tabcolsep}
 \begin{longtable}{p{\halfdeflength} | p{\halfdeflength}}
  \textbf{Dimension four:} Given a $(g;k_1,k_2,k_3)$ trisecting Morse $2$--function $f:X \to \R^2$, let $X_i = f^{-1}(A_i)$. Then $X=X_1 \cup X_2 \cup X_3$ is a $(g;k_1,k_2,k_3)$ trisection of $X$.
  &
  \textbf{Dimension three:} Given a genus $g$ Heegaard splitting Morse function $f:M \to \R$, let $M_1 = f^{-1}(R_\pi)$ and $M_2 = f^{-1}(R_0)$. Then $M=M_1 \cup M_2$ is a genus $g$ Heegaard splitting of $X$.
 \end{longtable}
\end{lemma}

As a hint for the proof, the only nonstandard part is to prove, in the four dimensional case, that each $X_i$ is diffeomorphic to a ``genus $k$'' $4$--dimensional $1$--handlebody $Z_k$. The best way to see this is to consider orthogonal projection from the sector $A_i \subset \R^2$, which is bounded by the rays $R_{2\pi (i-1)/3}$ and $R_{2\pi i/3}$, onto the intermediate ray $R_{\pi (2i-1)/3}$. After a suitable isotopy in $\R^2$ fixed along the two bounding rays, we may assume that the only places where the image of the singular locus is vertical with respect to this orthogonal projection are the midpoints of the indefinite folds which do not have cusps and the midpoint of the definite fold. Composing $f$ with this projection can then be seen to be a Morse function on $X_i$ with only critical points of index $0$ and $1$, with one index $0$ and $k$ index $1$ critical points.

Note that the distinct critical values condition (in dimension three) and double points avoiding cusps condition (in dimension four) are not strictly necessary to make the above lemma true, but we add them as a conceptual convenience.

The most basic example of a Heegaard splitting Morse function is the projection $(x_1,x_2,x_3,x_4) \mapsto x_1$ on $S^3 \subset \R^4$, giving the standard genus $0$ Heegaard splitting of $S^3$. Similary, projection $(x_1,x_2,x_3,x_4,x_5) \mapsto (x_1,x_2)$ is a trisecting Morse $2$--function on $S^4 \subset \R^5$, giving the $(0,0)$ trisection of $S^4$. The reader should verify these basic facts. It is also not too hard to see a Heegaard splitting Morse function on $S^1 \times S^2$ inducing a genus $1$ Heegaard splitting, and a trisecting Morse $2$--function on $S^1 \times S^3$ inducing a $(1,1)$ trisection. Beyond this, it is in fact not usually very straightforward to write down explicit Morse functions and Morse $2$--functions, let alone ones that induce the decompositions we desire. More frequently, we understand the decomposition first, from some other description of the manifold, and from this we can understand an appropriate Morse function or $2$--function.

The existence part of Theorem~\ref{T:ExUn} can, however, be proved by proving the existence of Heegaard splitting Morse functions and trisecting Morse $2$--functions. The former is standard, done by proving first the existence of Morse functions, then showing that one can cancel pairs of critical points until there is only one index $0$ and one index $3$ critical point, and finally showing that critical points can be rearranged so that their corresponding critical values increase with increasing index. A proof of the latter appears in~\cite{GayKirby} starting from a handle decomposition of the $4$--manifold, but can probably also be proved in a purely Morse $2$--function theoretic method, starting with the existence of Morse $2$--functions and then arguing that the critical locus of a Morse $2$--function can be cleaned up by a sequence of standard moves to become a trisecting Morse $2$--function. The work of Baykur and Saeki~\cite{BaykurSaeki} should provide enough tools to do this.

The uniqueness part of Theorem~\ref{T:ExUn} is proved in the three dimensional case using standard Cerf theory, where stabilization of the Heegaard splitting corresponds to adding a cancelling pair of index $1$ and $2$ critical points. See~\cite{Laudenbach} for a careful exposition of this proof. The four dimensional uniqueness proof in~\cite{GayKirby} unfortunately does not follow this parallel, i.e. does not use a Morse $2$--function version of Cerf theory, but is rather more ad hoc. For the sake of completeness it would be nice to see a Cerf theoretic proof, although it is not clear if the ultimate payoff would be worth the time.

One challenging but reachable example that the reader who likes working with explicit expressions might enjoy is to write down a trisecting Morse $2$--function on $\C P^2$. This can be done by suitably perturbing the following moment map:
\[[z_0:z_1:z_2] \mapsto \left(\frac{|z_1|^2}{|z_0|^2+|z_1|^2+|z_2|^2},\frac{|z_2|^2}{|z_0|^2+|z_1|^2+|z_2|^2}\right)\] 
The moment map itself is not a Morse $2$--function, but adding a generic perturbation term should make it Morse, and careful choice of this perturbation should make it trisecting. 

On the other hand, one can extract a trisection directly from this moment map without perturbing it to a Morse $2$--function. Let $x = |z_1|^2/(|z_0|^2+|z_1|^2+|z_2|^2)$ and $y=|z_2|^2/(|z_0|^2+|z_1|^2+|z_2|^2)$. Then the following ``tropical'' decomposition is in fact a $(1,0)$ trisection of $X=\C P^2$, and verifying this is also a good exercise:
\begin{align*} 
X_1 &= \{x \leq 1/4, y \leq 1/4\}\\
X_2 &= \{y \geq 1/4, y \geq x\}\\
X_3 &= \{x \geq 1/4, x \geq y\}
\end{align*}
This trisection is used in~\cite{LambertCole} to give a combinatorial proof of the Thom conjecture (that algebraic curves in $\C P^2$ minimize genus in their homology classes), and is also used in~\cite{LambertColeMeier} to produce minimal genus trisections of a large class of algebraic surfaces.

It is well known, from standard Morse theory, that stabilization of a Heegaard splitting can be realized via a homotopy of Morse functions supported in a ball, in which a cancelling pair of index $1$ and $2$ critical points is born near the Heegaard surface. Similarly, stabilization of a trisection can be realized via a homotopy of Morse $2$--functions supported in a ball; this is the ``introduction of an eye'' discussed in Section~5 of~\cite{GayKirby}. 

\section{The diagrammatic perspective} \label{S:Diagrams}

\begin{figure}
    \labellist
    \small\hair 2pt
    \pinlabel $\alpha_i$ [b] at 50 11
    \pinlabel $\alpha'_i$ [bl] at 83 80
    \pinlabel $\alpha_j$ [b] at 181 11 
    \pinlabel $a$ [b] at 114 50
    \endlabellist
 \centering
 \includegraphics{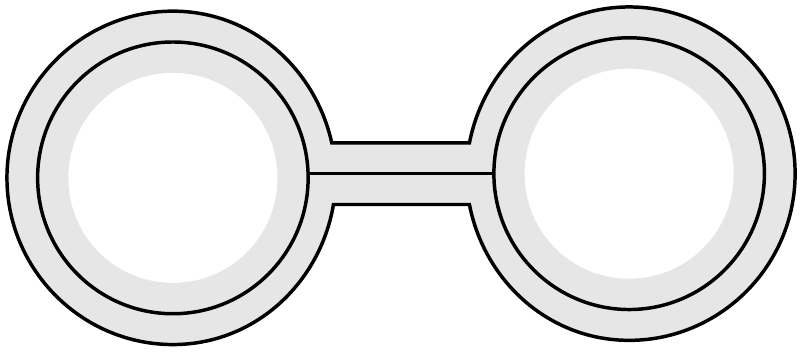}
 \caption{A handle slide.} \label{F:Slide}
\end{figure}
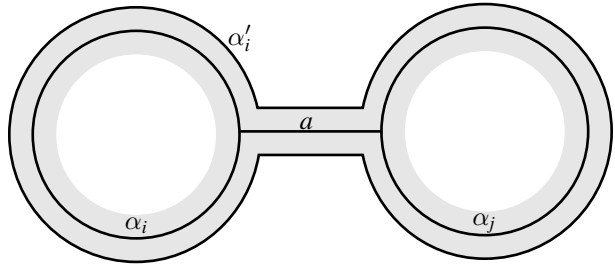

\begin{definition}
 Given a finite collection $\mathbf{\alpha}=\{\alpha_1,\ldots,\alpha_n\}$ of disjoint simple closed curves on an oriented surface $\Sigma$, given two of these curves $\alpha_i$ and $\alpha_j$ and an arc $a$ joining $\alpha_i$ and $\alpha_j$ and otherwise disjoint from the $\alpha$ curves, one can produce a new collection of curves $\mathbf{\alpha'}=\{\alpha'_1, \ldots, \alpha'_n\}$ by {\em sliding $\alpha_i$ over $\alpha_j$ along $a$} as follows: For $k \neq i$, $\alpha_k$ is unchanged, i.e. $\alpha'_k = \alpha_k$, while $\alpha_i'$ is the unique boundary component of a regular neighborhood of $\alpha_i \cup a \cup \alpha_j$ which is not isotopic to either $\alpha_i$ or $\alpha_j$ (see Figure~\ref{F:Slide}). Two collections of disjoint simple closed curves on the same surface are {\em slide equivalent} if one can be transformed to the other by a sequence of handle slides and isotopies. Two pairs of collections of simple closed curves $(\alpha,\beta)$ and $(\alpha',\beta')$ on the same surface are {\em slide equivalent} if $\alpha$ is slide equivalent to $\alpha'$ and $\beta$ is slide equivalent to $\beta'$. Two triples $(\Sigma,\alpha,\beta)$ and $(\Sigma',\alpha',\beta')$ are {\em slide diffeomorphic} if $(\alpha,\beta)$ is slide equivalent to some $(\alpha'',\beta'')$ such that $(\Sigma,\alpha'',\beta'')$ is orientation preserving diffeomorphic to $(\Sigma',\alpha',\beta')$.
\end{definition}

\begin{definition} \label{D:StdDiagrams}
In which we define Heegaard and trisection diagrams.
\setlength{\halfdeflength}{0.5\linewidth}
 \addtolength{\halfdeflength}{-2\tabcolsep}
 \begin{longtable}{p{\halfdeflength} | p{\halfdeflength}}
  \textbf{Dimension four:} A {\em $(g;k_1,k_2,k_3)$ trisection diagram} is a tuple $(\Sigma,\alpha,\beta,\gamma)$ where $\Sigma$ is a closed oriented surface of genus $g$ and the triples $(\Sigma,\alpha,\beta)$, $(\Sigma,\beta,\gamma)$ and $(\Sigma,\gamma,\delta)$ are each slide diffeomorphic to the standard Heegaard diagram $(\Sigma_g,\alpha^{g,k_i},\beta^{g,k_i})$ shown in Figure~\ref{F:StdDiagramsTri}. (Here $i=1$ for $(\alpha,\beta)$, $i=2$ for $(\beta,\gamma)$ and $i=3$ for $(\gamma,\alpha)$.)
  &
  \textbf{Dimension three:} A {\em genus $g$ Heegaard diagram} is a tuple $(\Sigma,\alpha,\beta)$ where $\Sigma$ is a closed oriented surface of genus $g$ and the pairs $(\Sigma,\alpha)$ and $(\Sigma,\beta)$ are both diffeomorphic to the standard pair $(\Sigma_g,\alpha^g)$ shown in Figure~\ref{F:StdDiagramsHS}.
 \end{longtable}
\end{definition}

\begin{figure}
  \begin{subfigure}[c]{\textwidth}
    \labellist
    \small\hair 2pt
    \pinlabel $k_i$ [t] at 74 8
    \pinlabel $g-k_i$ [t] at 218 10
    \endlabellist
    \centering
    \includegraphics{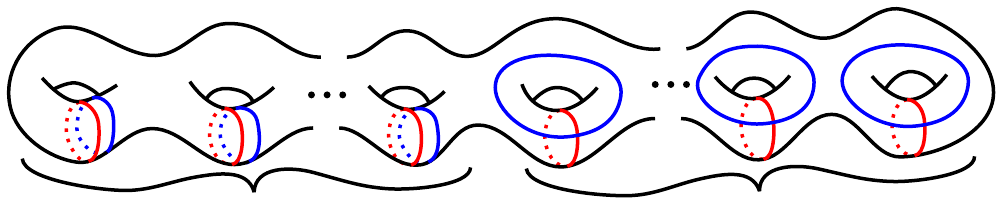}
    \caption{Standard Heegaard diagrams $(\Sigma_g,\alpha^{g,k_i},\beta^{g,k_i})$}
    \label{F:StdDiagramsTri}
  \end{subfigure}
  \begin{subfigure}[c]{\textwidth}
    \labellist
    \small\hair 2pt
    \pinlabel $g$ [t] at 73 6
    \endlabellist
    \centering
    \includegraphics{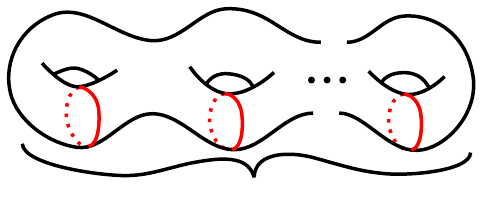}
    \caption{Standard genus $g$ cut system $(\Sigma_g,\alpha^g)$}
    \label{F:StdDiagramsHS}
  \end{subfigure}
  \caption{Standard diagrams needed for the definition of trisection and Heegaard diagrams. Red is $\alpha$ and blue is $\beta$.}
   \label{F:StdDiagrams}
\end{figure}

\pagebreak

Note that any system of curves $\alpha$ on a surface $\Sigma$ of genus $g$ such that $(\Sigma,\alpha) \cong (\Sigma_g, \alpha^g)$ is called a {\em cut system}, and cut systems are more standardly characterized by the property that they cut the ambient surface into a punctured sphere. We should think of $(\Sigma_g, \alpha^g)$ as the ``standard genus $g$ cut system''.

\begin{figure}
 \centering
 \includegraphics{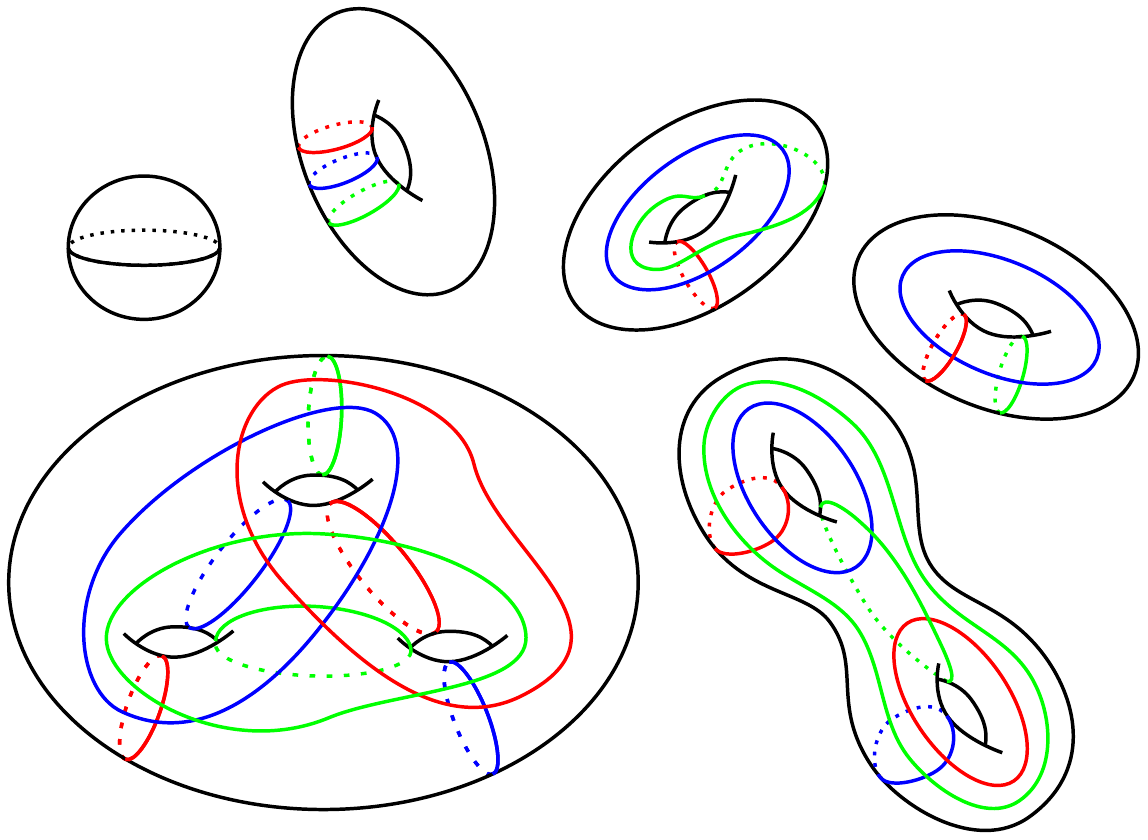}
 \caption{A selection of trisection diagrams. Red is $\alpha$, blue is $\beta$ and green is $\gamma$.} \label{F:TriDiagExamples}
\end{figure}

The reader should now verify that all the diagrams in Figure~\ref{F:TriDiagExamples} are in fact trisection diagrams. The genus $3$ example is the only one  in which handle slides are required to exhibit the standardness of pairs of colors. An instructive additional exercise is simply to try to draw new nontrivial diagrams from scratch; here ``nontrivial'' would first mean not diffeomorphic to connected sums of any of these examples, and at a second pass, not slide diffeomorphic to connected sums of these. Of special interest are diagrams which cannot be made ``simultaneously standard'' in the sense that no slides are needed to exhibit pairwise standardness. In thinking carefully about this problem one may discover the following fact: Simultaneously standard trisection diagrams are in fact connected sums of diagrams of genus $1$ and $2$. This is proved by thinking about the euler characteristic and number of boundary components of chains of curves each intersecting the next once.

\begin{lemma}
 In which we relate Heegaard and trisection diagrams to Heegaard splittings and trisections.
 \setlength{\halfdeflength}{0.5\linewidth}
 \addtolength{\halfdeflength}{-2\tabcolsep}
 \begin{longtable}{p{\halfdeflength} | p{\halfdeflength}}
  \textbf{Dimension four:} Given a trisection diagram $\Dd= (\Sigma,\alpha,\beta,\gamma)$ there is a $4$--manifold $X = X(\Dd)$ with trisection $\Tt(\Dd)= (X,X_1,X_2,X_3)$ such that $\Sigma = X_1 \cap X_2 \cap X_3$, oriented according to the conventions in Definition~\ref{D:HSandTri}, and such that the $\alpha$ curves bound embedded disks in $X_3 \cap X_1$, the $\beta$ curves in $X_1 \cap X_2$ and the $\gamma$ curves in $X_2 \cap X_3$.
  &
  \textbf{Dimension three:} Given a Heegaard diagram $\Dd = (\Sigma,\alpha,\beta)$ there is a $3$--manifold $M=M(\Dd)$ with Heegaard splitting $\Ss(\Dd) = (M,M_1,M_2)$ such that $\Sigma = M_1 \cap M_2$, oriented  according to the conventions in Definition~\ref{D:HSandTri}, and such that the $\alpha$ curves bound embedded disks in $M_1$ and the $\beta$ curves bound embedded disks in $M_2$.
  \\
  \textbullet\;  Any other trisected $4$--manifold satisfying these same properties with respect to the given diagram $\Dd$ is in fact orientation preserving diffeomorphic to $\Tt(\Dd)$.
  &
  \textbullet\; Any other Heegaard split $3$--manifold satisfying these same properties with respect to the given diagram $\Dd$ is in fact orientation preserving diffeomorphic to $\Ss(\Dd)$.
  \\
  \textbullet\; For every trisection $\Tt=(X,X_1,X_2,X_3)$ of a $4$--manifold $X$ there is a trisection diagram $\Dd$ such that $\Tt \cong \Tt(\Dd)$.
  &
  \textbullet\; For every Heegaard splitting $\Ss=(M,M_1,M_2)$ of a $3$--manifold $M$ there is a Heegaard diagram $\Dd$ such that $\Ss \cong \Ss(\Dd)$.
 \end{longtable}
\end{lemma}

\pagebreak

The proof is again left as an exercise. In dimension three, the key point is that any collection of curves on $\Sigma$ diffeomorphic to $(\Sigma_g,\alpha^g)$ (a cut system) uniquely determines a handlebody filling of $\Sigma$ in which these curves bound disks. In dimension four we need this fact, the fact that the standard Heegaard diagram $(\Sigma^g,\alpha^{g,k},\beta^{g,k})$ used in the definition of trisection diagram is a diagram for $\#^k (S^1 \times S^2)$, and the fact~\cite{LaudPoen} that any self diffeomorphism of $\#^k (S^1 \times S^2)$ extends to a self diffeomorphism of $\#^k (S^1 \times B^3)$.

Now we consider how the uniqueness statement for Heegaard splittings of $3$--manifolds and trisections of $4$--manifolds translates into the world of diagrams. 

\begin{definition}
 In which we define certain standard diagrams used in the stabilization process.
\setlength{\halfdeflength}{0.5\linewidth}
 \addtolength{\halfdeflength}{-2\tabcolsep}
 \begin{longtable}{p{\halfdeflength} | p{\halfdeflength}}
  \textbf{Dimension four:} The standard $(1;1,0,0)$ trisection diagram for $S^4$ is the diagram $\Dd^*_1 = (T^2,\mu,\mu,\lambda)$ shown in Figure~\ref{F:StableTri}. Cyclically permuting the curve systems gives the standard $(1;0,1,0)$ diagram $\Dd^*_2$  and the standard $(1;0,0,1)$ diagram $\Dd^*_3$.
  &
  \textbf{Dimension three:} The standard genus $1$ Heegaard diagram for $S^3$ is the diagram $\Dd^* = (T^2,\mu,\lambda)$ shown on the right in Figure~\ref{F:StableHS}.
 \end{longtable}
\end{definition}
\begin{figure}
  \begin{subfigure}[c]{.5\textwidth}
    \centering
    \includegraphics{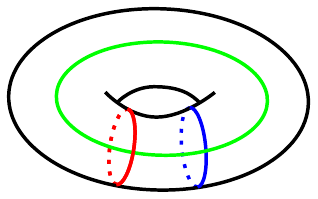}
    \caption{Stabilizing trisection diagram $\Dd^*_1$}
    \label{F:StableTri}
  \end{subfigure}
  \begin{subfigure}[c]{.5\textwidth}
    \centering
    \includegraphics{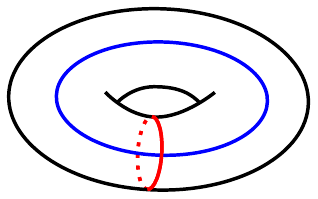}
    \caption{Stabilizing Heegaard diagram $\Dd^*$}
    \label{F:StableHS}
  \end{subfigure}
  \caption{Standard stabilization diagrams}
   \label{F:StableDiagrams}
\end{figure}
The implication in the above definitions is of course that $\Ss(\Dd_1)$ is a Heegaard splitting of $S^3$ and that $\Tt(\Dd_{(1;1,0,0)})$ is a trisection of $S^4$ (and the same for the other two standard trisection diagrams). The reader should verify these facts, and also verify the following:

\begin{proposition} \label{P:StabilizingDiagrams}
 In which we relate stabilization of splittings and trisections to diagrams.
\setlength{\halfdeflength}{0.5\linewidth}
 \addtolength{\halfdeflength}{-2\tabcolsep}
 \begin{longtable}{p{\halfdeflength} | p{\halfdeflength}}
  \textbf{Dimension four:} Given a trisection diagram $\Dd$ with associated trisected $4$--manifold $\Tt = \Tt(\Dd)$, let $\Tt'$ be the result of an $i$--stabilization of $\Tt$. Then $\Tt' \cong \Tt(\Dd \# \Dd^*_i)$.
  &
  \textbf{Dimension three:} Given a Heegaard diagram $\Dd$ with associated Heegaard split $3$--manifold $\Ss = \Ss(\Dd)$, let $\Ss'$ be the result of stabilizing $\Ss$. Then $\Ss' \cong \Ss(\Dd \# \Dd_1)$.
  \\
  \textbullet\; Given two trisection diagrams $\Dd$ and $\Dd'$, with 
  \[
  \Tt(\Dd) = (X,X_1,X_2,X_3)
  \]
  and
  \[
  \Tt(\Dd') = (X',X'_1,X'_2,X'_3),
  \]
we have that $X \cong X'$ if and only if, for some $k_1, k_2, k_3$ and $k'_1,k'_2,k'_3$, the following two trisection diagrams are slide diffeomorphic:
\[\Dd \# (\#^{k_1} \Dd^*_1) \# (\#^{k_2} \Dd^*_2) \# (\#^{k_3} \Dd^*_3)\] 
and 
\[\Dd' \# (\#^{k'_1} \Dd^*_1) \# (\#^{k'_2} \Dd^*_2) \# (\#^{k'_3} \Dd^*_3)\]
  &
  \textbullet\; Given two Heegaard diagrams $\Dd$ and $\Dd'$, with \[\Ss(\Dd) = (M,M_1,M_2)\] 
  and 
  \[\Ss(\Dd') = (M',M'_1,M'_2),\] 
  we have that $M \cong M'$ if and only if, for some $k$ and $k'$, the following two Heegaard diagrams are slide diffeomorphic:
  \[\Dd \# (\#^k \Dd^*)\] 
  and 
  \[\Dd' \# (\#^{k'} \Dd^*)\]
 \end{longtable}
 
\end{proposition}

Now we discuss trisection diagrams in relation to trisecting Morse $2$--functions; the phenomena discussed here are unique to the $4$--dimensional setting and do not have obvious $3$--dimensional analogues. Recall the notation from Section~\ref{S:Morse}, in particular the ``trisection'' of $\R^2$ as $\R^2 = A_1 \cup A_2 \cup A_3$, the labelling of rays $R_\theta$ by angle $-\theta$ to the positive $x$--axis, and the identification of each ray with $[0,\infty)$. Fix a closed $4$--manifold $X$ with a trisecting Morse $2$--function $f: X \to \R^2$, and consider the induced trisection $X = X_1 \cup X_2 \cup X_3$, where $X_i = f^{-1}(A_i)$, as in Lemma~\ref{L:Fcn2Decomp}. We can read off a trisection diagram from the function $f$ and a (generic) choice of gradient-like vector field over each ray $R_{2\pi i/3}$, since this data gives us descending manifolds for each of the index $2$ critical points in each handlebody, i.e. a handle decomposition of each handlebody, with the attaching ``spheres'' in the central surface $\Sigma = X_1 \cap X_2 \cap X_3$ being a collection of simple closed curves $\alpha$, $\beta$ or $\gamma$.

What we would like to emphasize here is that there is more information available in a Morse $2$--function than just in the trisection diagram. In fact, if we choose a smoothly varying family of gradient-like vector fields over the rays $R_\theta$, i.e. on each $f^{-1}(R_\theta)$ we choose a gradient-like vector field $V_\theta$ for the Morse functions $f_\theta : R_\theta \to [0,\infty)$, smooth in $\theta$, then we can look at the descending manifolds for the index $2$ critical points of $f_\theta$ in $\Sigma = f_\theta^{-1}(0)$. There we will see a ``moving family'' of cut systems, mostly moving by isotopy but, at isolated times, experiencing discrete moves. More precisely, from a trisecting Morse $2$--function on a closed $4$--manifold $X$, we can first arrange that all the cusps in each sector appear at the same $\theta$ value, and then choosing one representative $\theta$--value during each $\theta$ interval when only isotopies occur, we can produce an {\em augmented trisection diagram} 
\[ (\Sigma, \alpha^1, \ldots, \alpha^a, \beta^1, \ldots, \beta^b, \gamma^1, \ldots, \gamma^c) \]
satisfying the following properties:
\begin{itemize}
 \item For any indices $i$, $j$ and $k$, $(\Sigma, \alpha^i, \beta^j, \gamma^k)$ is a trisection diagram for $X$.
 \item The cut system $\alpha^{i+1}$ is obtained from $\alpha^i$ by a single handle slide, and similary for the $\beta$'s and $\gamma$'s.
 \item The Heegaard diagrams $(\Sigma,\alpha^a,\beta^1)$, $(\Sigma,\beta^b,\gamma^1)$ and $(\Sigma,\gamma^c,\alpha^1)$ are each {\em diffeomorphic} (not just slide diffeomorphic) to the standard diagram $(\Sigma_g, \alpha^{g,k_i}, \beta^{g,k_i})$.
\end{itemize}
(If we think of cut systems as {\em ordered lists} of simple closed curves, ordered by the relative heights of the corresponding critical points, rather than just as {\em sets} of simple closed curves, then we should also include transposition of two adjacent curves in the list as a move, and this would correspond to one critical point rising above another.)

In what sense does this augmented trisection diagram contain more information than an ordinary trisection diagram? The main point is that rather than simply asserting that each pair $(\Sigma,\alpha,\beta)$, $(\Sigma,\beta,\gamma)$ and $(\Sigma,\gamma,\alpha)$ is slide diffeomorphic to a standard Heegaard diagram, with the augmented diagram we now know exactly {\em how} to slide handles to get each pair to be standard. This in turn means that, rather than simply appealing to Laudenbach and Po\'enaru to assert that we can fill in each sector with $\natural^k S^1 \times B^3$, i.e. with $3$--handles and a $4$--handle, and that any way of filling in is as good as any other way, we actually explicitly see the attaching maps for the $3$--handles. Also, the minimum length $a+b+c$ of an augmented trisection diagram for a given trisection is a measure of the complexity of the trisection, which should be of interest and is closely related to certain complexity measures coming from simplicial complexes associated to curve systems on surfaces; see~\cite{KirbyThompson} for example.

\section{$4$--manifolds with boundary}

One advantage to thinking of trisections from a Morse $2$--function perspective is that this gives us the most natural definition of a trisection of a $4$--manifold with boundary. 
\begin{definition}
In Definition~\ref{D:HSTfunctions}, we assumed that the manifolds we were working with were closed. Now suppose they have nonempty boundary instead.
\setlength{\halfdeflength}{0.5\linewidth}
\addtolength{\halfdeflength}{-2\tabcolsep}
\setlength{\wholedeflength}{2\halfdeflength}
\addtolength{\wholedeflength}{2\tabcolsep}
 \begin{longtable}{p{\halfdeflength} | p{\halfdeflength}}
  \textbf{Dimension four:} A {\em $(g,k)$ trisecting} Morse $2$--function $f$ on a $4$--manifold $X$ with nonempty connected boundary is a Morse $2$--function $f: X \to \R^2$ such that:
  &
  \textbf{Dimension three:} A {\em genus $g$ Heegard splitting} Morse function $f$ on a $3$--manifold $M$ with nonempty connected boundary is a Morse function $f:M \to \R$ such that:
  \\
  \textbullet\; $f(X) = D^2$.
  &
  \textbullet\; $f(M) = [-1,1]$.
  \\
  \textbullet\; For all points $p \in S^1 = \partial D^2$, $f^{-1}(p)$ is a compact surface with boundary, contained in $\partial X$, and in fact the restriction of $f$ to $f^{-1}(S^1)$ is a compact surface bundle over $S^1$.
  &
  \textbullet\; $f^{-1}(1)$ and $f^{-1}(-1)$ are diffeomorphic compact surfaces with boundary, contained in $\partial M$.
  \\
  \textbullet\; The closure of the complement of $f^{-1}(S^1)$ in $\partial M$ is diffeomorphic to $B \times D^2$, for $B$ a collection of circles and with $f$ being projection onto the $D^2$ factor. For each $p \in S^1$, $B \times \{p\}$ is the boundary of $f^{-1}(p)$.
  &
  \textbullet\; The closure of the complement of $f^{-1}(\{-1,1\})$ in $\partial M$ is diffeomorphic to $B \times [-1,1]$, for $B$ a collection of circles and with $f$ being projection onto the $[-1,1]$ factor. $B \times \{-1\}$ is the boundary of $f^{-1}(-1)$ and $B \times \{1\}$ is the boundary of $f^{-1}(1)$.
  \\
  \textbullet\; $\mathbf{0}=(0,0)$ is a regular value of $f$, and thus $f^{-1}(\mathbf{0}) = \Sigma$ is a compact surface with boundary, where $\partial \Sigma = B \times \{\mathbf{0}\} \subset B \times D^2$.
  &
  \textbullet\; $0$ is a regular value of $f$, and thus $f^{-1}(0) = \Sigma$ is a compact surface with boundary, where $\partial \Sigma = B \times \{0\} \subset B \times [-1,1]$.
  \\
  \textbullet\; On each of the three rays $R_0$, $R_{2\pi/3}$ and $R_{4\pi/3}$, $f$ has only index $2$ critical points, all of which have distinct critical values which lie in the interiors of the rays (with the same number of critical points on each ray).
  &
  \textbullet\; On each of the two rays $R_0$ and $R_\pi$, $f$ has only index $2$ critical points, all of which have distinct critical values which lie in the interiors of the rays (with the same number of critical points on each ray).
  \\
  \cmidrule{2-2}
  \multicolumn{2}{p{\wholedeflength}}{\textbullet\; Over each of the three sectors $A_1$, $A_2$ and $A_3$, each component of the singular locus of $f$ is an arc from one bounding ray of $A_i$ to the next, with at most one cusp per component. All folds are indefinite folds. Furthermore, in $\R^2$ each of these components is transverse to each ray $R_\theta$ except at cusps, which are tangent to the rays, and $f$ restricted to the singular locus is an immersion with cusps and double points avoiding the cusps.}
 \end{longtable}
\end{definition}

From this we can give the following definition:
\begin{definition}
Let $X$, resp. $M$, be a compact $4$--manifold, resp. $3$--manifold, with nonempty connected boundary.
\setlength{\halfdeflength}{0.5\linewidth}
 \addtolength{\halfdeflength}{-2\tabcolsep}
 \begin{longtable}{p{\halfdeflength} | p{\halfdeflength}}
  \textbf{Dimension four:} A {\em relative trisection} of $X$ is a decomposition $X = X_1 \cup X_2 \cup X_3$ for which there exists a trisecting Morse $2$--function $f: M \to \R^2$ with $X_i = f^{-1}(A_i)$.
  &
  \textbf{Dimension three:} A {\em sutured Heegaard splitting} of $M$ is a decomposition $M=M_1 \cup M_2$ for which there exists a Heegaard splitting Morse function $f: M \to \R$ with $M_1 = f^{-1}(R_\pi)$ and $M_2 = f^{-1}(R_0)$.
 \end{longtable}
\end{definition}

These are not the standard definitions, but we feel that the Morse theoretic perspective better conveys the central idea. The standard definitions can be recovered with some observations/exercises:
\begin{itemize}
 \item Starting in \textbf{dimension three}, the induced structure on $\partial M$ is a decomposition $\partial M = -R_- \cup ([-1,1] \times \partial R_-) \cup R_+$, where $R_-$ and $R_+$ are diffeomorphic oriented compact surfaces with boundary. A $3$--manifold with such a structure on its boundary is a balanced sutured $3$--manifold.
 \item The two pieces $M_1$ and $M_2$ can each be viewed as either {\em sutured compression bodies} from the central surface $\Sigma = f^{-1}(0)$ to $R_\pm$ or as handlebodies, where $\partial M_1 = \Sigma \cup ([-1,0] \times \partial \Sigma) \cup -R_-$ and $\partial M_2 = -\Sigma \cup ([0,1] \times \partial \Sigma) \cup R_+$.
 \item The $3$--dimensional part of definition~\ref{D:Stab}, stabilization of Heegaard splittings, makes sense in this relative setting, assuming the stabilizing arc lies entirely in the interior of $M$. The assertion that the result is again a sutured Heegaard splitting, using the Morse theoretic definition of splitting given above, requires seeing that stabilization is achieved by perturbing the Morse function to introduce a cancelling $1$--$2$ critical point pair.  
 \item Moving to \textbf{dimension four}, the induced structure on $\partial X$ is an open book decomposition, namely a decomposition into a surface bundle over $S^1$ ($E \subset \partial X$ with $f: E \to S^1$) and a disjoint union of solid tori $B \times D^2$, such that the boundary of each fiber $f^{-1}(\theta)$ of the surface bundle (each ``page'') is the link $B \times \{\theta\}$ in $B \times D^2$. These pages are traditionally extended by adding on the annuli $B \times R_\theta$, to get Seifert surfaces for the link $B \times \mathbf{0}$, the ``binding'' of the open book.
 \item Each pairwise intersection $X_i \cap X_j$ is a sutured compression body from the central surface $\Sigma = f^{-1}(\mathbf{0})$ to the page $f^{-1}(e^{2\pi i/3})$ (identifying $\R^2$ with $\C$).
 \item Each piece $X_i$ is a $4$--dimensional $1$--handlebody, but it's boundary comes with a decomposition into three pieces: $\partial X_i = (X_i \cap X_{i-1}) \cup (X_i \cap X_{i+1}) \cup (X_i \cap \partial X)$. The first two pieces are the above mentioned sutured compression bodies, and the third part $X_i \cap \partial X$ is one third of the open book decomposition of $\partial X$.
 \item The {\em internal portion} of $\partial X_i$, i.e. the closure of $\partial X_i \setminus \partial X$, comes equipped with a sutured Heegaard splitting, i.e $(X_i \cap X_{i-1}) \cup (X_i \cap X_{i+1})$.
 \item In fact each such $X_i$ is diffeomorphic to $C_i \times [-1,1]$ for some sutured compression body $C_i$ from some surface $\Sigma'$ to the page $f^{-1}(e^{2\pi i/3})$, with the internal portion of $\partial X_i$ being diffeomorphic to $(C_i \times \{-1\}) \cup (\Sigma' \times [-1,1]) \cup (C_i \times \{1\})$.
 \item Note that the preceding item also gives a sutured Heegaard splitting of the internal portion of $\partial X_i$, namely as the union of $(C_i \times \{-1\}) \cup (\Sigma' \times [-1,0])$ and $(C_i\times \{1\}) \cup (\Sigma' \times [0,1])$. The previous Heegaard splitting $(X_i \cap X_{i-1}) \cup (X_i \cap X_{i+1})$ is a stabilization of this Heegaard splitting.
 \item The $4$--dimensional version of stabilization in Definition~\ref{D:Stab} also now makes sense, again assuming the stabilizing lies entirely in the interior of $X$. Seeing that the result is again a relative trisection according to our definition of relative trisections as coming from trisecting Morse $2$--functions requires seeing that stabilization can be achieved by a perturbation of the Morse $2$--function. The stabilization takes place in a neighborhood of a boundary parallel arc in some $X_{i-1} \cap X_{i+1}$. The perturbation of the Morse $2$--function takes place in a neighborhood of that arc and its boundary parallelizing disk, and literally pulls the arc from lying over $A_{i-1} \cap A_{i+1}$ back to lying over $A_i$. This is discussed in detail in~\cite{GayKirby}.
\end{itemize}

The fundamental existence and uniqueness result from the closed case still holds in this relative setting, provided we work relative to fixed boundary data:

\begin{theorem}[Existence and Uniqueness] \label{T:ExUnRel}
 The above decompositions exist and are unique up to stabilization relative to fixed boundary data. More precisely:
 \setlength{\halfdeflength}{0.5\linewidth}
 \addtolength{\halfdeflength}{-2\tabcolsep}
 \begin{longtable}{p{\halfdeflength} | p{\halfdeflength}}
  \textbf{Dimension four:} Given any open book decomposition on the boundary $\partial X$ of a compact connected oriented $4$--manifold with nonempty connected boundary, there exists a relative trisection of $X$ inducing this open book on $\partial X$. Any two trisections of the same $4$--manifold inducing the same open book on the boundary become isotopic after some number of stabilizations~\cite{GayKirby}.
  &
  \textbf{Dimension three:} Given any balanced sutured decomposition of the boundary of a compact connected oriented $3$--manifold $M$ with nonempty connected boundary, there exists a sutured Heegaard splitting on $M$ inducing the sutured structure on $\partial M$. Any two sutured Heegaard splittings of the same $3$--manifold inducing the same boundary data become isotopic after some number of stabilizations.
 \end{longtable}
\end{theorem}

\begin{theorem}[Gluing]
  These relative decompositions are especially useful because they can be glued together when the boundary data agree.
 \setlength{\halfdeflength}{0.5\linewidth}
 \addtolength{\halfdeflength}{-2\tabcolsep}
 \begin{longtable}{p{\halfdeflength} | p{\halfdeflength}}
  \textbf{Dimension four:} Given relatively trisected $4$--manifolds $X = X_1 \cup X_2 \cup X_3$ and $X' = X'_1 \cup X'_2 \cup X'_3$ and an orientation reversing diffeomorphism $\phi : \partial X \to \partial X'$ respecting the induced open book decompositions, then the following decomposition of the closed $4$--manifold $\tilde{X} = X \cup_\phi X'$ is a trisection~\cite{Castro}:
  \[ \tilde{X} = (X_1 \cup_\phi X'_1) \cup (X_2 \cup_\phi X'_2) \cup (X_3 \cup_\phi X'_3) \]
  &
  \textbf{Dimension three:} Given $3$--manifolds with sutured Heegaard splittings $M=M_1 \cup M_2$ and $M'=M'_1 \cup M'_2$ and an orientation reversing diffeomorphism $\phi: \partial M \to \partial M'$ respecting the induced sutured decompositions, then the following decomposition of the closed $3$--manifold $\tilde{M} = M \cup_\phi M'$ is a Heegaard splitting:
  \[ \tilde{M}= (M_1 \cup_\phi M'_1) \cup (M_2 \cup_\phi M'_2) \]
 \end{longtable}
\end{theorem}

The reader should prove the $3$--dimensional statement as an exercise. The $4$--dimensional statement takes more work.

An important example of the boundary data one might consider comes, in both cases, when studying a knot complement. 

\begin{figure}
 \centering
 \includegraphics{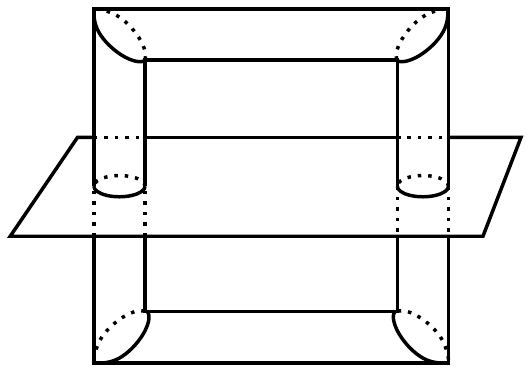}
 \caption{A knot in $1$--bridge position giving a sutured Heegaard splitting.}\label{F:SuturedKnotComplement}
\end{figure}

In dimension three, a classical knot $K$ in $S^3$ gives rise to its exterior $E(K) = S^3 \setminus \nu(K)$, a $3$--manifold with boundary parametrized as $S^1 \times S^1$, where the first $S^1$ factor is the meridian, i.e. the boundary of $D^2$ in $\nu(K) \cong D^2 \times S^1$. Then (see Figure~\ref{F:SuturedKnotComplement}) identifying the second $S^1$ factor as $\partial ([-1,1] \times [-1,1])$, we can decompose $\partial E(K)$ as $-R_- \cup (\partial R_i \times [-1,1]) \cup R_+$ where $R_- = S^1 \times [-1,1] \times \{-1\}$, $R_+ = S^1 \times [-1,1] \times \{1\}$ and $\partial R_- \times [-1,1] = S^1 \times \{-1,1\} \times [-1,1]$. A sutured Heegaard splitting of $E(K)$ with respect to these sutures is precisely the restriction of an ordinary Heegaard splitting of $S^3$ to $E(K)$ when $E(K)$ is in $1$--bridge position with respect to this splitting (equivalently, when $K$ is represented by a doubly pointed Heegaard diagram). Also, it is not important here that the knot $K$ is in $S^3$, the same construction works in any closed $3$--manifold. But in $S^3$ this is the standard construction used to apply sutured Floer homology to knot complements~\cite{Juhasz}.

In dimension four, a knotted sphere $K: S^2 \into S^4$ has exterior $E(K)$ with $\partial E(K)$ parametrized as $S^1 \times S^2$, again with the first $S^1$ factor being the meridian, i.e. the boundary of $D^2$ in $\nu(K) \cong D^2 \times S^2$. Now if we identify $S^2$ with 
$\partial ([-1,1] \times D^2)$, we get a natural open book decomposition on $\partial E(K)$, where the surface bundle part is the annulus bundle $S^1 \times [-1,1] \times S^1$, projecting onto the second $S^1$ factor, and the neighborhood of the binding is the union of two solid tori $S^1 \times \{-1,1\} \times D^2$. We can use the existence theorem above to conclude that this open book extends to a relative trisection of the exterior, but in fact~\cite{GayMeier} shows that such a relative trisection is actually the restriction to $E(K)$ of an ordinary trisection of $S^4$ when the sphere $K$ is in ``$1$--bridge trisection position'' with respect to this trisection, and~\cite{MeierZupanGenBridge} shows how to put any sphere in such a position. As above, being in $S^4$ is not essential.

It remains to discuss relative trisections from the diagrammatic perspective. 
\begin{definition} \label{D:RelDiagrams}
In which we define sutured Heegaard diagrams and relative trisection diagrams.
\setlength{\halfdeflength}{0.5\linewidth}
 \addtolength{\halfdeflength}{-2\tabcolsep}
 \begin{longtable}{p{\halfdeflength} | p{\halfdeflength}}
  \textbf{Dimension four:} A relative trisection diagram is a tuple $(\Sigma,\alpha,\beta,\gamma)$ where $\Sigma$ is a compact oriented surface with nonempty boundary and the triples $(\Sigma,\alpha,\beta)$, $(\Sigma,\beta,\gamma)$ and $(\Sigma,\gamma,\delta)$ are each slide diffeomorphic to a standard sutured Heegaard diagram as shown in Figure~\ref{F:StdDiagramsRelTri}.
  &
  \textbf{Dimension three:} A sutured Heegaard diagram is a tuple $(\Sigma,\alpha,\beta)$ where $\Sigma$ is a compact oriented surface with nonempty boundary and the pairs $(\Sigma,\alpha)$ and $(\Sigma,\beta)$ are both diffeomorphic to a standard pair as shown on the right in Figure~\ref{F:StdDiagramsRelHS}.
 \end{longtable}
\end{definition}
\begin{figure}
  \begin{subfigure}[c]{\textwidth}
    \centering
    \includegraphics{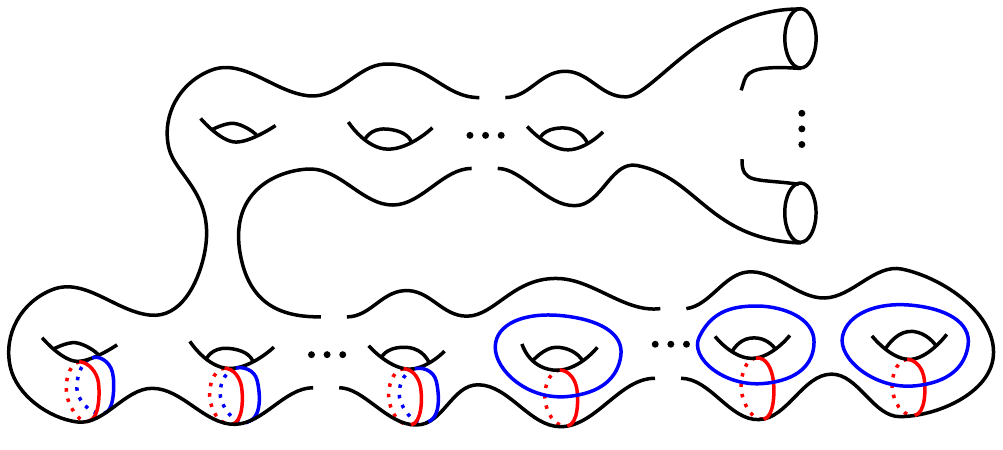}
    \caption{Standard sutured Heegaard diagrams for definition of relative trisection diagram.}
    \label{F:StdDiagramsRelTri}
  \end{subfigure}
  \begin{subfigure}[c]{\textwidth}
    \centering
    \includegraphics{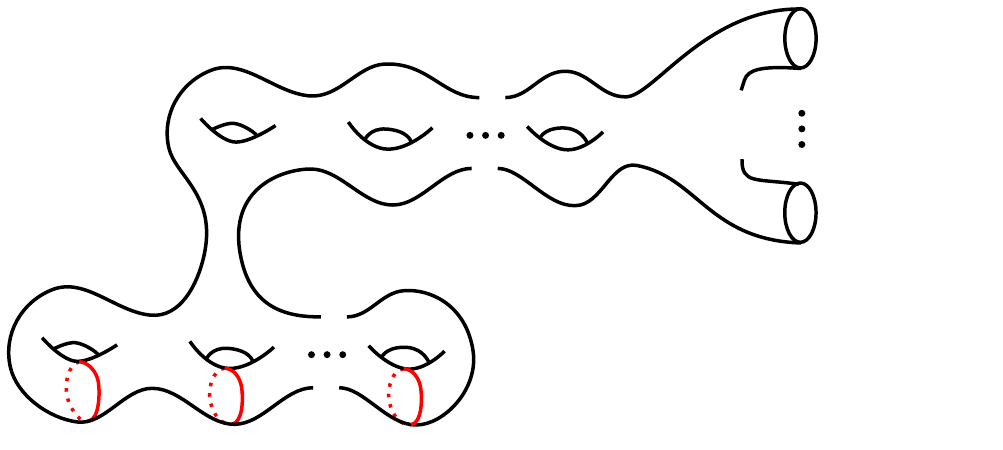}
    \caption{Standard curve system for definition of sutured Heegaard diagram.}
    \label{F:StdDiagramsRelHS}
  \end{subfigure}
  \caption{Standard diagrams needed for the definition of relative trisection and sutured Heegaard diagrams.}
   \label{F:StdRelDiagrams}
\end{figure}
Note that we have dropped reference to the genus and other parameters describing the exact standard diagrams used for these definitions since the naming of the parameters at this point seems not to be helpful. Other references take care to name the genus, number of boundary components, number of curves, and so forth.

Parts of the following result can be proved as basic exercises in both dimensions, but in dimension four the heart of the result is perhaps nontrivial and is proved in~\cite{CastroGayPinzon}, to which the reader is referred.
\begin{proposition}
 In which we relate sutured Heegaard diagrams and relative trisection diagrams to sutured Heegaard splittings and relative trisections.
 \setlength{\halfdeflength}{0.5\linewidth}
 \addtolength{\halfdeflength}{-2\tabcolsep}
 \begin{longtable}{p{\halfdeflength} | p{\halfdeflength}}
  \textbf{Dimension four:} Given a relative trisection diagram $\Dd= (\Sigma,\alpha,\beta,\gamma)$ there is a compact $4$--manifold $X = X(\Dd)$ with nonempty connected boundary with relative trisection $\Tt(\Dd)= (X,X_1,X_2,X_3)$ such that $\Sigma = X_1 \cap X_2 \cap X_3$, oriented according to the conventions in Definition~\ref{D:HSandTri}, and such that the $\alpha$ curves bound embedded disks in $X_3 \cap X_1$, the $\beta$ curves in $X_1 \cap X_2$ and the $\gamma$ curves in $X_2 \cap X_3$.
  &
  \textbf{Dimension three:} Given a sutured Heegaard diagram $\Dd = (\Sigma,\alpha,\beta)$ there is a sutured $3$--manifold $M=M(\Dd)$ with sutured Heegaard splitting $\Ss(\Dd) = (M,M_1,M_2)$ such that $\Sigma = M_1 \cap M_2$, oriented according to the conventions in Definition~\ref{D:HSandTri}, and such that the $\alpha$ curves bound embedded disks in $M_1$ and the $\beta$ curves bound embedded disks in $M_2$.
  \\
  \textbullet\;  Any other relatively trisected $4$--manifold satisfying these same properties with respect to the given diagram $\Dd$ is in fact orientation preserving diffeomorphic to $\Tt(\Dd)$.
  &
  \textbullet\; Any other $3$--manifold with a sutured Heegaard splitting satisfying these same properties with respect to the given diagram $\Dd$ is in fact orientation preserving diffeomorphic to $\Ss(\Dd)$.
  \\
  \textbullet\; For every relative trisection $\Tt=(X,X_1,X_2,X_3)$ of a $4$--manifold $X$ there is a relative trisection diagram $\Dd$ such that $\Tt \cong \Tt(\Dd)$.
  &
  \textbullet\; For every sutured Heegaard splitting $\Ss=(M,M_1,M_2)$ of a $3$--manifold $M$ there is a sutured Heegaard diagram $\Dd$ such that $\Ss \cong \Ss(\Dd)$.
 \end{longtable}
\end{proposition}

Coming full circle to the Morse theoretic perspective, the last assertion in the result above, that sutured Heegaard splittings and relative trisections come from diagrams, can be shown by seeing that a Morse function or Morse $2$--function inducing the given decomposition yields, via the appropriate gradient-like vector fields, descending manifolds for the index $2$ critical points that intersect the central surface in precisely the curves of the diagram.
 
\section{Surfaces in $4$--manifolds}

Meier and Zupan in~\cite{MeierZupanBridge} introduced the notion of bridge trisections of surfaces embedded in $S^4$, as the natural trisected generalization of bridge splittings of knots in $S^3$, and further generalized this to embedded surfaces in arbitrary $4$--manifolds. Following our theme, we introduce these decompositions from a Morse theoretic point of view. To do this we begin with something easier than Definition~\ref{D:MorseMorse2}:

\begin{definition} \label{D:MorseMorse2LowD}
 In which we define Morse functions and Morse $2$--functions in the limited context of dimensions one and two.
 \setlength{\halfdeflength}{0.5\linewidth}
 \addtolength{\halfdeflength}{-2\tabcolsep}
 \setlength{\wholedeflength}{2\halfdeflength}
 \addtolength{\wholedeflength}{2\tabcolsep}
 \begin{longtable}{p{\halfdeflength} | p{\halfdeflength}}
  \textbf{Dimension two:} A {\em Morse $2$--function} on a $2$--manifold $S$ is a smooth function $f: S \to \R^2$ which, at every point $p \in S$, has one of the following three forms with respect to appropriate local coordinates $(t,x)$ near $p$ and $(u,v)$ near $f(p)$:
  &
  \textbf{Dimension one:} A {\em Morse function} on a $1$--manifold $K$ is a smooth function $f: K \to \R$ which, at every point $p \in M$, has one of the following two forms with respect to appropriate local coordinates $x$ near $p$ and $u$ near $f(p)$:
  \\
  \textbullet\; $(t,x) \mapsto (u=t,v=x)$; here $p$ is called a {\em regular point}.
  &
  \textbullet\; $x \mapsto u=x$; here $p$ is called a {\em regular point}.
  \\
  \textbullet\; $(t,x) \mapsto (u=t, v = \pm x^2)$; here $p$ is called a {\em fold point}; in this dimension all folds are definite.
  &
  \textbullet\; $x \mapsto u = \pm x^2$; here $p$ is called a {\em critical point}, of index $0$ if $u=x^2$ and index $1$ if $u=-x^2$.
  \\\cmidrule{2-2}
  \multicolumn{2}{p{\wholedeflength}}{\textbullet\; $(t,x) \mapsto (u=t, v= x^3 - t x)$; here $p$ is called a {\em cusp point}.}
 \end{longtable}
\end{definition}

Here are the parallel basic facts to check (again assume that the domain of $f$ is closed):
\begin{itemize}
 \item In both cases the inverse image of a regular value is an even number of points.
 \item In both cases the singular locus, the set of all critical points, is a closed codimension one submanifold, i.e. a finite collection of points in dimension one and a finite collection of embedded circles in dimension two.
 \item In dimension two, the cusp points form a finite collection of points on the singular locus.
 \item Via a small perturbation one may assume that the critical points of a Morse function have distinct critical values
 \item Letting $Z$ be the singular locus of a Morse $2$--function $f$, via a small perturbation one may assume that $f|_Z$ is an immersion with semicubical cusps, with at worst double point self intersections, none of which occur at cusps.
 \item If $f:S \to \R^2$ is a Morse $2$--function and $A$ is an arc in $\R^2$ avoiding the cusps and transverse to the image of the singular locus, then $K = f^{-1}(A)$ is a $1$--manifold in $X$, with $\partial K = f^{-1}(\partial A)$. 
 \item Furthermore, if we identify $A$ with an interval in $\R$ via some embedding $A \hookrightarrow \R$ then $f|_M :M \to A$ is a Morse function with critical points where $A$ crosses folds. Reversing the orientation of $A$ changes the indices of these critical points, with index $0$ becoming index $1$ and vice versa. 
 \item Crossing a definite fold in the index $0$ direction adds a new pair of points to the fiber while crossing in the index $1$ direction removes such a pair.
\end{itemize}

\begin{definition} \label{D:HSTfunctionsLowD}
In which we define bridge splitting Morse functions and bridge trisecting Morse $2$--functions. (Recall the $R_\theta$ and $A_i$ notation introduced earlier for rays and sectors in $\R^2$.)
\setlength{\halfdeflength}{0.5\linewidth}
\addtolength{\halfdeflength}{-2\tabcolsep}
\setlength{\wholedeflength}{2\halfdeflength}
\addtolength{\wholedeflength}{2\tabcolsep}
\begin{longtable}{p{\halfdeflength} | p{\halfdeflength}}
  \textbf{Dimension two:} A {\em bridge trisecting} Morse $2$--function $f$ on a surface $S$ is a Morse $2$--function $f: S \to \R^2$ such that:
  &
  \textbf{Dimension one:} A {\em bridge splitting} Morse function $f$ on a $1$--manifold $K$ is a Morse function $f:K \to \R$ such that:
  \\
  \textbullet\; $\mathbf{0}=(0,0)$ is a regular value of $f$.
  &
  \textbullet\; $0$ is a regular value of $f$.
  \\
  \textbullet\; On each of the three rays $R_0$, $R_{2\pi/3}$ and $R_{4\pi/3}$, $f$ has only index $1$ critical points, all of which have distinct critical values.
  &
  \textbullet\; On each of the two rays $R_0$ and $R_\pi$, $f$ has only index $1$ critical points, all of which have distinct critical values.
  \\\cmidrule{2-2}
  \multicolumn{2}{p{\wholedeflength}}{\textbullet\; Over each of the three sectors $A_1$, $A_2$ and $A_3$, the singular locus of $f$ consists of arcs from one bounding ray of $A_i$ to the next with at most one cusp on each arc. Furthermore, in $\R^2$ each of these components is transverse to each ray $R_\theta$ except at cusps, which are tangent to the rays, and $f$ restricted to the singular locus is an immersion with cusps and double points avoiding the cusps.}
 \end{longtable}
\end{definition}

A bridge splitting function on a $1$--manifold $K$ decomposes $K$ into $K_1 \cup K_2$, where each $K_i$ is a collection or arcs. A bridge trisecting function on a surface $S$ decomposes $S$ into $S_1 \cup S_2 \cup S_3$, where each $S_i$ is a disjoint union of disks, each $S_i \cap S_j$ is a disjoint union of arcs, and $S_1 \cap S_2 \cap S_3$ is an even number of points.
\begin{definition} \label{D:HSTfunctionsPairs}
In which we define bridge splitting Morse functions and bridge trisecting Morse $2$--functions on pairs.
\setlength{\halfdeflength}{0.5\linewidth}
 \addtolength{\halfdeflength}{-2\tabcolsep}
 \begin{longtable}{p{\halfdeflength} | p{\halfdeflength}}
  \textbf{Dimensions two and four:} A {\em bridge trisecting} Morse $2$--function $f$ on a surface $S$ embedded in a $4$--manifold $X$ is a trisecting Morse $2$--function $f:X \to \R^2$ such that $f|_S$ is a bridge trisecting Morse $2$--function on $S$.
  &
  \textbf{Dimensions one and three:} A {\em bridge splitting} Morse function $f$ on a $1$--manifold $K$ embedded in a $3$--manifold $M$ is a Heegaard splitting Morse function $f:M \to \R$ such that $f|_K$ is a bridge splitting Morse function on $K$.
  \\
  \textbullet\; A {\em generalized bridge trisection} of a surface $S$ embedded in a $4$--manifold $X$ is a decomposition $(X,S) = (X_1,S_1) \cup (X_2,S_2) \cup (X_3,S_3)$ coming from a bridge trisecting function in the the sense that $X_i = f^{-1}(A_i)$ and $S_i = S \cap X_i$.
  &
  \textbullet\; A {\em generalized bridge splitting} of a knot or link $K$ in a $3$--manifold $M$ is a decomposition $(M,K) = (M_1,K_1) \cup (M_2,K_2)$ coming from a bridge splitting function $f$, in the sense that $M_1 = f^{-1}(R_\pi)$, $M_2 = f^{-1}(R_0)$ and $K_i = K \cap M_i$.
  \\
  \textbullet\; A {\em bridge trisection} is a generalized bridge trisection of an embedded surface in $S^4$ for which the underlying trisection of $S^4$ has genus $0$.
  &
  \textbullet\; A {\em bridge splitting} is a generalized bridge splitting of a knot or link in $S^3$ for which the underlying Heegaard splitting has genus $0$.
 \end{longtable}
\end{definition}
The reader show now check that the following conditions follow; these are usually taken as the standard definitions of bridge splitting and bridge trisection:
\begin{enumerate}
 \item In dimensions one and three, the arcs making up $K_i$ are properly embedded and simultaneously boundary parallel in the handlebody $M_i$.
 \item In dimensions two and four, the disks making up $S_i$ are properly embedded and simultaneously boundary parallel in the $4$--dimensional $1$--handlebody $X_i$.
 \item The intersection $S_i \cap S_j$ is a collection of arcs properly embedded and simultaneously boundary parallel in the handlebody $X_i \cap X_j$.
\end{enumerate}

Meier and Zupan~\cite{MeierZupanBridge} showed that every surface in $S^4$ can be isotoped so as to be bridge trisected by the standard genus $0$ trisecting Morse $2$--function, and later~\cite{MeierZupanGenBridge} showed how to do this in arbitrary $4$--manifolds with respect to arbitrary trisections. The analogous statement for knots and links in $3$--manifolds is standard. There are also uniqueness statements up to stabilization moves, but we will not discuss those here.

There are actually several interesting ways to think about (generalized) bridge trisections diagrammatically. Since there is quite a lot to say, we describe these vaguely and give references for details. Honest bridge trisections, and bridge splittings, are described by tangles in $3$--balls, so these can simply be drawn as standard tangle diagrams. These are the diagrams discussed in~\cite{MeierZupanBridge}, and are called {\em triplane diagrams}. Trivial (boundary parallel) tangles can also be described as half-plat closures of braids, and thus bridge trisections can also be described by braids; this perspective is important in Saltz's work~\cite{Saltz} on Khovanov-style invariants of surfaces in $S^4$. In more general $4$--manifolds, one needs to record the trisection of the $4$--manifold as well as the surface, and this can be done either through multi-pointed diagrams or by ``shadow diagrams'', in which each arc in each tangle is described by it's shadow on the trisecting surface. Shadow diagrams are used in~\cite{LambertColeMeier}, while multi-pointed diagrams are discussed in~\cite{GayMeier}.

As a final remark, in the discussion in this section we have assumed that the ambient manifolds and submanifolds are closed; the fully relative case, in which either or both may have boundary, is obviously more subtle but can be understood with care. The details have been worked out by Meier~\cite{MeierInPrep}.

%
%
%
\bibliographystyle{plain}
%

%

\bibliography{WinterBraidsNotes}

\end{document}